\documentclass[11pt]{amsart}
\usepackage{amsmath,amsthm,amsfonts,amssymb,latexsym,epsfig}

\newtheorem{theorem}{Theorem}[section]

\newtheorem{lemma}{Lemma}[section]

\newtheorem{cor}{Corollary}[section]

\numberwithin{equation}{section}

\theoremstyle{definition}

\theoremstyle{remark}

\begin{document}
\title{On weighted mean matrices whose $l^p$ norms are determined on decreasing sequences}
\author{Peng Gao}
\address{Division of Mathematical Sciences, School of Physical and Mathematical Sciences,
Nanyang Technological University, 637371 Singapore}
\email{penggao@ntu.edu.sg}
\date{October 6, 2008.}


\begin{abstract}
  We give a condition on weighted mean matrices so that their $l^p$ norms are determined on decreasing sequences when the condition is satisfied.
  We apply our result to give a proof of a conjecture of Bennett and discuss some related results.
\end{abstract}

\maketitle
\section{Introduction}
\label{sec 1} \setcounter{equation}{0}

  Suppose throughout that $p\neq 0, \frac{1}{p}+\frac{1}{q}=1$.
  For $p \geq 1$, let $l^p$ be the Banach space of all complex sequences ${\bf a}=(a_n)_{n \geq 1}$ with norm
\begin{equation*}
   ||{\bf a}||_p: =(\sum_{n=1}^{\infty}|a_n|^p)^{1/p} < \infty.
\end{equation*}
  The celebrated
   Hardy's inequality (\cite[Theorem 326]{HLP}) asserts that for $p>1$,
\begin{equation}
\label{eq:1} \sum^{\infty}_{n=1}\Big{|}\frac {1}{n}
\sum^n_{k=1}a_k\Big{|}^p \leq \Big (\frac
{p}{p-1} \Big )^p\sum^\infty_{n=1}|a_n|^p.
\end{equation}
   Hardy's inequality can be regarded as a special case of the
   following inequality:
\begin{equation}
\label{01}
   \Big | \Big |C \cdot {\bf a}\Big | \Big |^p_p =\sum^{\infty}_{n=1} \Big{|}\sum^{\infty}_{k=1}c_{n,k}a_k
    \Big{|}^p \leq U_p \sum^{\infty}_{n=1}|a_n|^p,
\end{equation}
   in which $C=(c_{n,k})$ and the parameter $p>1$ are assumed
   fixed, and the estimate is to hold for all complex
   sequences ${\bf a} \in l^p$. The $l^{p}$ operator norm of $C$ is
   then defined as
\begin{equation*}
\label{02}
    ||C||_{p,p}=\sup_{||{\bf a}||_p = 1}\Big | \Big |C \cdot {\bf a}\Big | \Big |_p.
\end{equation*}
   It follows that inequality \eqref{01} holds for any ${\bf a} \in l^p$ when $U^{1/p}_p \geq ||C||_{p,p}$ and fails to hold for some ${\bf a} \in l^p$
   when $U^{1/p}_p <||C||_{p,p}$.
    Hardy's inequality thus asserts that the Ces\'aro matrix
    operator $C$, given by $c_{n,k}=1/n , k\leq n$ and $0$
    otherwise, is bounded on {\it $l^p$} and has norm $\leq
    p/(p-1)$. (The norm is in fact $p/(p-1)$.)

    We say a matrix $A=(a_{n,k})$ is a lower triangular matrix if $a_{n,k}=0$ for $n<k$ and a lower triangular matrix $A$ is a summability matrix if
    $a_{n,k} \geq 0$ and
    $\sum^n_{k=1}a_{n,k}=1$. We say a summability matrix $A$ is a weighted
    mean matrix if its entries satisfy:
\begin{equation}
\label{021}
    a_{n,k}=\lambda_k/\Lambda_n,  ~~ 1 \leq k \leq
    n; \hspace{0.1in} \Lambda_n=\sum^n_{i=1}\lambda_i, \hspace{0.1in} \lambda_i \geq 0, \hspace{0.1in} \lambda_1>0.
\end{equation}
     We shall also say that a weighted mean matrix $A$ is generated
     by $\{ \lambda_n \}^{\infty}_{n=1}$ (resp. $\{ \lambda_n
     \}^{N}_{n=1}$) when $A$ is an infinite weighted mean matrix
     (resp. finite $N \times N$ weighted mean matrix) whose
     entries are given by \eqref{021}.

    Hardy's inequality \eqref{eq:1} motivates one to
    determine the $l^{p}$ operator norm of an arbitrary summability or weighted mean matrix $A$.
    In the weighted mean matrix case, as the diagonal entries $\{ \lambda_n/\Lambda_n \}$ uniquely determine one such
    a matrix, one certainly expects to obtain a bound for its norm using only the diagonal terms. In \cite{G5}, the author proved the following result:
\begin{theorem}
\label{thm03}
    Let $1<p<\infty$ be fixed. Let $A$ be a weighted mean matrix generated
     by $\{ \lambda_n \}^{\infty}_{n=1}$ . If for any integer $n \geq 1$, there exists a positive constant
    $0<L<p$ such that
\begin{equation}
\label{024}
    \frac {\Lambda_{n+1}}{\lambda_{n+1}} \leq \frac
    {\Lambda_n}{\lambda_n}  \Big (1- \frac
    {L\lambda_n}{p\Lambda_n} \Big )^{1-p}+\frac {L}{p}~~,
\end{equation}
    then
    $||A||_{p,p} \leq p/(p-L)$.
\end{theorem}

   It is easy to see that the above result implies the following well-known result of Cartlidge \cite{Car} (see also \cite[p. 416, Theorem C]{B1}):
\begin{theorem}
\label{thm02}
    Let $1<p<\infty$ be fixed. Let $A$ be a weighted mean matrix generated
     by $\{ \lambda_n \}^{\infty}_{n=1}$ . If
\begin{equation}
\label{022}
    L=\sup_n\Big(\frac {\Lambda_{n+1}}{\lambda_{n+1}}-\frac
    {\Lambda_n}{\lambda_n}\Big) < p ~~,
\end{equation}
    then
    $||A||_{p,p} \leq p/(p-L)$.
\end{theorem}

   The above result of Cartlidge is often very handy to apply for determining
   $l^p$ norms of certain weighted mean matrices, when combined with a result of Cass and Kratz \cite{CF}, which says that for
   a weighted mean matrix $A$ generated by $\{ \lambda_n
   \}^{\infty}_{n=1}$, with the $\lambda_n$'s generated by a positive logarithmico-exponential
   function (for details, see \cite{G}) and satisfying $\lim_{n \rightarrow
   \infty}\Lambda_n/(n\lambda_n)=L<p$, then $||A||_{p,p} \geq
   p/(p-L)$. As an example, we note the following two
     inequalities were claimed to hold (with no proofs supplied) by Bennett ( \cite[p. 40-41]{B4}; see also \cite[p. 407]{B5}):
\begin{eqnarray}
\label{7}
   \sum^{\infty}_{n=1}\Big{|}\frac
1{n^{\alpha}}\sum^n_{i=1}(i^{\alpha}-(i-1)^{\alpha})a_i\Big{|}^p &
\leq & \Big( \frac {\alpha p}{\alpha p-1} \Big )^p\sum^{\infty}_{n=1}|a_n|^p, \\
\label{8}
   \sum^{\infty}_{n=1}\Big{|}\frac
1{\sum^n_{i=1}i^{\alpha-1}}\sum^n_{i=1}i^{\alpha-1}a_i\Big{|}^p &
\leq & \Big(\frac {\alpha p}{\alpha p-1} \Big
)^p\sum^{\infty}_{n=1}|a_n|^p,
\end{eqnarray}
     whenever $p>1, \alpha p >1$. We note here the constant $(\alpha p /(\alpha p-1))^p$ is best possible by the result of Cass and Kratz (or see \cite{Be1}).

   Straightforward applications of Theorem
\ref{thm02} allow the author \cite{G} to prove inequalities
\eqref{7} for $p>1, \alpha \geq 1$ and
     \eqref{8} for $p>1, \alpha \geq 2$ or $0< \alpha \leq 1, \alpha p >1$. The same result was obtained for \eqref{8}
     by Bennett himself
\cite{Be1} independently and his proof also relies on Cartlidge's
result. Using a different approach, Bennett was able to prove
\eqref{7} for the full range of $\alpha$ (see
     \cite[Theorem 1]{Be1} with $\beta=1$ there). Using the result of Theorem \ref{thm03}, the author
\cite{G5} has shown that inequality \eqref{8} holds for $p \geq 2,
1< \alpha <2$ (in fact, as pointed out in \cite{G5}, for fixed
$1<p<2$, one can also prove \eqref{8} for some cases of
$1<\alpha<2$).

   We note here that by a change of variables $a_k \rightarrow a^{1/p}_k$ in \eqref{eq:1} and on letting $p \rightarrow +\infty$, one obtains the
   following well-known Carleman's inequality  \cite{Carlman}, which asserts that for convergent infinite series $\sum a_n$ with non-negative terms,
   one has
\begin{equation*}
   \sum^\infty_{n=1}(\prod^n_{k=1}a_k)^{\frac 1{n}}
\leq e\sum^\infty_{n=1}a_n,
\end{equation*}
   with the constant $e$ being best possible.

   It is then natural to study the following weighted version of Carleman's inequality:
\begin{equation}
\label{1}
   \sum^N_{n=1}\Big( \prod^n_{k=1}a^{\lambda_k/\Lambda_n}_k \Big )
\leq E_N\sum^N_{n=1}a_n,
\end{equation}
  where the notations are as in \eqref{021} and $N \geq 1$ is an integer or $N=\infty$.
  The task here is to determine the best constant $E_N$ so that inequality \eqref{1} holds
  for any (convergent when $N=\infty$) series $\sum a_n$ with non-negative terms. Note that
  \eqref{1} can be regarded as the $p \rightarrow +\infty$ case of the following
  inequality (once again by a change of variables):
\begin{equation}
\label{3.1}
    \sum^{N}_{n=1} \Big{(}\sum^{n}_{k=1}\frac {\lambda_k}{\Lambda_n}a_k
    \Big{)}^p \leq U_{p, N} \sum^{N}_{n=1} a_n^p,
\end{equation}
   where $U_{p, N}$ is a positive constant, $a_n \geq 0$ and
   $\lambda_n, \Lambda_n$'s are given as in \eqref{021}.

    Note that Cartlidge's result (Theorem \ref{thm02}) implies that when \eqref{022} is satisfied, then for any ${\bf a} \in l^p$,
  inequality \eqref{3.1} holds for any $N$ with $U_{p, N}=(p/(p-L))^p$.
  Similar to our discussions above, by a change of variables $a_k \rightarrow a^{1/p}_k$ in \eqref{3.1} and on letting $p \rightarrow +\infty$,
  one obtains inequality \eqref{1} with $E_N=e^{L}$ as long as \eqref{022} is satisfied with $p$ replaced by $+\infty$ there.

  In connection to \eqref{8}, Bennett \cite[p. 829]{Be1} further conjectured that inequality \eqref{1} holds
   for $\lambda_k=k^{\alpha}$ for $\alpha > -1$ with $E_{\infty}=e^{1/(\alpha+1)}$. As the cases $-1 < \alpha \leq 0$ or $\alpha \geq 1$
   follow directly from the known cases of inequalities \eqref{8} upon changes of
   variables $\alpha \rightarrow \alpha+1, a_k \rightarrow a^{1/p}_k$ and on letting $p \rightarrow +\infty$,
   the only nontrivial cases are when $0< \alpha <1$. As these
   cases are the limits of the corresponding $l^p$ cases
   and the author \cite{G5} has shown \eqref{8} hold for $p \geq 2, 1 < \alpha <2$ using Theorem \ref{thm03},
   it follows that Bennett's conjecture is true.

   Motivated by the study of inequalities \eqref{7}-\eqref{8}, we
   seek for extra inputs that may lead to a resolution of the remaining
   case of \eqref{8} for $1<p<2, 1 < \alpha <2$. For this, we note
   the following natural question related to the $l^p$ norms of any matrix asked by Bennett
  \cite[Problem 7.23]{B5}: When is the norm of a matrix
  determined by its action on decreasing sequences? In other
  words, when do we have
\begin{equation}
\label{1.7}
  ||C||_{p,p} = \sup \Big \{||C \cdot {\bf a}||_p: ||{\bf
a}||_p = 1 \text{ and  ${\bf a}$ decreasing} \Big \} ?
\end{equation}

   For weighted mean matrices, it is known that \cite[p. 422]{B1}
   that sequences $\bf{a}$, with $a_n/\lambda^{1/(p-1)}_n$
   decreasing in $n$, are sufficient to determine the norm.
   Note that this certainly implies \eqref{1.7} when the
   $\lambda_n$'s are decreasing. A slightly generalization of this
   later case is given in the following lemma:
\begin{lemma}\cite[Lemma 2.4]{CLO}
\label{lem3}
   Let $p>1$ and $C = (c_{n,k})_{n,k \geq 1}$ be an arbitrary lower triangular matrix.
If $c_{n,k} \geq  c_{n,k+1} \geq 0$ for all $n \geq 1, 1 \leq k
\leq n-1$, then \eqref{1.7} holds.
\end{lemma}

   We refer the reader to the articles \cite{CHS} and \cite{CSW} for more recent
   developments in this area. It is our goal in this paper to give a condition on weighted mean
matrices in Section \ref{sec 8'} so that \eqref{1.7} will hold. As
an application, we will give another proof of the above mentioned
Bennett's conjecture.

   We note that Cartlidge's result (Theorem \ref{thm02}) only allows one to prove
\eqref{7} with some restrictions on the
   $\alpha$'s, as in \cite{G},
   leaving alone the cases $1/p< \alpha \leq 1$. However, for these cases, Lemma \ref{lem3} implies that \eqref{1.7} holds
   for the corresponding matrices. This extra information can be used to give a proof of these cases and in fact we shall
    prove a more general result in Section \ref{sec 4}.

    In \cite{G5}, the author has shown that several approaches in the literature
   concerning the $l^p$ norms of weighted mean matrices are
   equivalent. In Section \ref{sec 8}, we will consider another approach to the $l^p$ norms of weighted mean
    matrices, namely the Schur's test. We will show that Schur's
   test is equivalent to the other approaches mentioned in \cite{G5} and we shall point out how Bennett's proof of \eqref{7} can be
rewritten using Schur's test. We shall also apply Schur's test to
give extensions of \eqref{7} which in turn allows us to view both
inequalities \eqref{7} and \eqref{8} as special cases of a family
of inequalities.

\section{On the validity of \eqref{1.7} for weighted mean matrices}
\label{sec 8'} \setcounter{equation}{0}

   In this section, we want to first present a result regarding the validity of \eqref{1.7} for weighted mean matrices.
   Since one can often reduce the questions of finding the norms of infinite weighted mean matrices to that of finite ones, we
   consider only finite weighted mean matrices here. Thus instead
   of \eqref{01}, we consider \eqref{3.1} instead and we have
\begin{theorem}
\label{thm2}
  Let $p>1$ be fixed and let $N \geq 1$ be a fixed integer and $A$ a weighted mean matrix
  generated by $\{ \lambda_n \}^N_{n=1}$. Suppose that \eqref{3.1} is satisfied for some positive constant $U_{p, N}$.
  If for any $1 \leq k \leq N-1$, the following condition
\begin{equation}
\label{1.8}
  \frac {1}{\Lambda_k} \geq U_{p, N} \Big(\frac {1}{\lambda_k}-\frac {1}{\lambda_{k+1}} \Big )
\end{equation}
  is satisfied, then \eqref{1.7} holds for $C=A$ in this case.
\end{theorem}
\begin{proof}
    Since our matrix $A$ is of finite dimension, it is easy to see that in this case we have
\begin{equation*}
    \mu^{1/p}_{p,N}:=||A||_{p,p}=\max_{||{\bf a}||_p = 1}\Big | \Big |A \cdot {\bf a}\Big | \Big |_p.
\end{equation*}
 Thus without loss of generality, we may assume that the maximum is reached at some ${\bf a}$ with $||{\bf a}||_p=1$. It is shown in \cite{G5} that
 in this case we have $a_n > 0$ for all $1 \leq n \leq N$ and on setting
\begin{equation*}
   A_n=\sum^n_{k=1}\frac {\lambda_ka_k}{\Lambda_n},
\end{equation*}
  we also have
\begin{equation}
\label{2.1'}
  \mu_{p,N} (\frac {a^{p-1}_k}{\lambda_k}-\frac {a^{p-1}_{k+1}}{\lambda_{k+1}})=\frac {A^{p-1}_k}{\Lambda_k}, \hspace{0.1in} 1 \leq k \leq N-1;
  \hspace{0.1in}  \mu_{p,N} \frac {a^{p-1}_N}{\lambda_N}=\frac {A^{p-1}_N}{\Lambda_N}; \hspace{0.1in} \sum^N_{n=1}a^p_n=1.
\end{equation}

  We now show by induction on $k$ that if \eqref{1.8} is satisfied, then the sequence ${\bf a}$ satisfying \eqref{2.1'} must be decreasing.
  First, it is easy to see that $a_1 \geq a_2$ using the relation $k=1$ in \eqref{2.1'} and noting that $A_1=a_1$ and $0<\mu_{p,N} \leq U_{p,N}$ by assumption.
  It now follows by induction that $A_k \geq a_k$ for $k \geq 1$ and that $a_k \geq a_{k+1}$ now follows from the $k$-th relation in \eqref{2.1'} and
  this establishes our assertion.
\end{proof}

  We note here that one sees from \eqref{2.1'} that
   that sequence $\bf{a}$ with $a_n/\lambda^{1/(p-1)}_n$
   decreasing in $n$, are sufficient to determine the norm, this is mentioned in Section \ref{sec
   1}.

   Now to apply Theorem \ref{thm2}, one needs to find some
   constant $U_{p,N}$ so that \eqref{3.1} holds. This is not a problem
   in many cases, as one can apply Theorem \ref{thm03} or Theorem
   \ref{thm02}. For example, if we use Theorem \ref{thm02}, then
   we can deduce the following result from Theorem \ref{thm2}:
\begin{cor}
\label{cor2}
  Let $p>1$ be fixed and let $N \geq 1$ be a fixed integer and $A$ a weighted mean matrix
  generated by $\{ \lambda_n \}^N_{n=1}$. Suppose that \eqref{022} is satisfied and for any $1 \leq k \leq N-1$,
  we have
\begin{equation}
\label{3.4}
  \Big (1-\frac {L}{p} \Big )^p \geq  \Lambda_k \Big(\frac {1}{\lambda_k}-\frac {1}{\lambda_{k+1}} \Big
  ),
\end{equation}
  then \eqref{1.7} holds for $C=A$ in this case.
\end{cor}

  We note that the left-hand side expression of \eqref{3.4} is an increasing
  function of $p$ for fixed $L$. Thus if $L<1$, then upon taking
  $p=1$, we see that \eqref{1.7} holds for any $p>1$ as long as
\begin{equation}
\label{4.3}
     \inf_n\Big(\frac {\Lambda_{n+1}}{\lambda_{n+1}}-\frac
    {\Lambda_n}{\lambda_n}\Big) \geq L.
\end{equation}
   One should compare the above with \eqref{022}. Interestingly enough, \eqref{4.3} tells us that if the condition \eqref{022} fails
   in the worst possible way (so that \eqref{4.3} holds), then Cartlidge's result (Theorem \ref{thm02}) does not help in determining the
   norm but then we can have an extra input by knowing that in this case \eqref{1.7} holds, provided that we know the norm is bounded by $p/(p-L)$. In particular, we point
   out here that if inequalities \eqref{8} were true for $p>1, 1<\alpha
   <2$ (note that it is shown in \cite{G5} that this is the case when $p \geq 2$), then Theorem \ref{thm2} implies that one
   may focus on decreasing sequences when trying to prove \eqref{8}, since in this case \eqref{4.3} holds with $\lambda_k=k^{\alpha-1}$ and $L=1/\alpha$
   (see \cite[Theorem 6]{Be1}). Of course one is not able to apply \eqref{1.8} using the constant $U_{p,N}=(\alpha p /(\alpha p-1))^p$ for the unknown cases of
   \eqref{8}.
    However,
   for the case of $p$ being large, one may hope to find a coarse bound
   $U_{p,N}$ so that \eqref{8} hold with the constant $(\alpha p /(\alpha
   p-1))^p$ replaced by $U_{p,N}$ and \eqref{1.8} is also satisfied and
   hopefully the extra information (that one may focus on
   decreasing sequences) will allow one to give a proof of
   \eqref{8} for the cases $1<\alpha <2$ and $p$
   large. We shall not worry about finding such a coarse bound
   here but we will show later in this section that the $p \rightarrow +\infty$ case (corresponding to the conjecture of Bennett mentioned
   in Section \ref{sec 1}) follows from this approach.

   By looking at the case $k=1$ of \eqref{2.1'}, we see that the
   case $k=1$ of \eqref{1.8} with $U_{p,N}$ replaced by $\mu_{p,N}$ is
   a necessary condition for $a_2 \geq a_1$. When $A=(a_{i,j})$ is an
   infinite weighted mean matrix, then we denote $A_N=(a_{i,j})_{1 \leq i, j \leq N}$ and let $\mu_{p,N} =||A_N||^{p}_{p,p}$ and note
   that we have $\mu_{p,N-1} \leq \mu_{p,N}$ for $N \geq 2$ (one sets $a_{N}=0$ in \eqref{3.1} to see
   this), thus the sequence $\{\mu_{p,N} \}^{\infty}_{N=1}$ is increasing
   and thus we have $\mu_{p,N} \rightarrow
   ||A||^{p}_{p,p}$ as $N \rightarrow +\infty$, which allows us to deduce
   immediately the following
\begin{cor}
\label{cor3}
  Let $p>1$ be fixed and $A$ a weighted mean matrix
  generated by $\{ \lambda_n \}^N_{n=1}$. A necessary
  condition for \eqref{1.7} to hold for $C=A$ is
\begin{equation*}
  \frac {1}{\lambda_1} \geq ||A||^{p}_{p,p} \Big(\frac {1}{\lambda_1}-\frac {1}{\lambda_{2}} \Big
  ).
\end{equation*}
   If moreover, the sequence
   $\{ \Lambda_n/\lambda_n \}^{\infty}_{n=1}$
   is convex, then the above condition is also sufficient.
\end{cor}

   We note here by a result of Bennett \cite[Theorem 2]{Be1}, we know that the
   sequence $\{ \Lambda_n/\lambda_n \}^{\infty}_{n=1}$ is convex when
   $\lambda_n=n^{\alpha}$ for $\alpha \geq 1$ or $\alpha \leq 0$ and is concave for $0 \leq \alpha \leq
   1$.

   We now consider two analogues of Theorem \ref{thm2} here.
   First we note that we have a similar result concerning inequality \eqref{1},
   namely,
\begin{theorem}
\label{thm3}
  Let $N \geq 1$ be a fixed integer and suppose that $E_{N}$ is the best possible constant to make \eqref{1} hold.
  If for any $1 \leq k \leq N-1$, inequality \eqref{1.8} is
  satisfied with $U_{p, N}$ replaced by $E'_N$ for some constant $E'_N \geq E_N$
  there, then to prove \eqref{1}, it suffices to establish it for decreasing
  sequences.
\end{theorem}

   Next, we note that one can also study inequality \eqref{3.1}
   when $p<0$ and one often expects to get result analogue to the
   case $p>0$. To be precise, we consider the following inequality
   for $a_n \geq 0$ and $p<0$,
\begin{equation}
\label{2.4}
    \sum^{N}_{n=1} \Big{(}\sum^{n}_{k=1}\frac {\lambda_k}{\Lambda_n}a^{1/p}_k
    \Big{)}^p \leq U_{p, N} \sum^{N}_{n=1} a_n.
\end{equation}
   Here we define the value of the left-hand side expression above
   to be $0$ when one or more of the $a_n$'s is zero. This makes
   the left-hand side expression above a continuous function on
   the compact set $\{ \sum^N_{n=1}a_n=1 | a_n \geq 0 \}$ and
   therefore we have $U_{p,N} < \infty$. From now on, for a weighted mean matrix $A$ generated by $\{ \lambda_n \}^N_{n=1}$
   ($N$ finite or infinite) and a fixed $p<0$, we shall
   denote $||A||^p_{p,p}$ for the supreme of the left-hand side
   expression of \eqref{2.4}, over the set $\{ \sum^N_{n=1}a_n=1 | a_n \geq 0
   \}$. We now have the following analogue of
   Cartlidge's result for $p<0$, which can be easily established by following the proof for the case $p>1$ given in \cite{G1}
   by noting that the case $n=1$ of \eqref{022} implies $L \geq
   0$.
\begin{theorem}
\label{thm1'}
  Let  $p<0$ be fixed and $A$ a weighted mean matrix generated by $\{ \lambda_n \}^N_{n=1}$.
  Then
\begin{equation*}
    \sum^{\infty}_{n=1}A^p_n
\geq \frac {p}{p-L} \sum^\infty_{n=1}a_nA^{p-1}_n,
\end{equation*}
   where $L$ is given as in \eqref{022}. In particular, inequality \eqref{2.4} holds with $||A||^p_{p, p} \leq (p/(p-L))^p$.
\end{theorem}

   Now, analogue to Theorem \ref{thm2}, we have
\begin{theorem}
\label{thm4}
  Let $p<0$ be fixed and $N \geq 1$ a fixed integer and $A$ a weighted mean matrix
  generated by $\{ \lambda_n \}^N_{n=1}$ and suppose that \eqref{2.4}
  holds for some constant $U_{p, N}$. If for any $1 \leq k \leq N-1$, inequality \eqref{1.8} is
  satisfied with $U_{p, N}$, then $||A||^p_{p,p}$ is determined on an increasing
  sequence.
\end{theorem}

   Now, we want to see what can be
   said about the $l^p$ norm of a given matrix, taking into the account that \eqref{1.7} holds for such a matrix. One strategy is to
   find a matrix whose $l^p$ norm (or an upper bound of it) is known,
   say by Cartlidge's result. Then one can make a comparison of
   the two matrices, thanks to the following result:
\begin{lemma}{\cite[Lemma 2.1]{B5}}
\label{lem2}
    Let ${\bf u}, {\bf v}$ be $n$-tuples
 with non-negative entries with $n \geq 1$ and
\begin{equation*}
   \sum_{i=1}^k u_i \leq \sum_{i=1}^k v_i, \hspace{0.1in} 1 \leq k \leq n-1; \hspace{0.1in} \sum_{i=1}^n u_i = \sum_{i=1}^n v_i.
\end{equation*}
    then
\begin{equation*}
   \sum_{i=1}^n u_i a_i \leq \sum_{i=1}^n v_i a_i,
\end{equation*}
   for any decreasing $n$-tuple ${\bf a}$ and the above inequality reverses when ${\bf a}$ is increasing.
\end{lemma}

   We note that the above lemma is given in \cite[Lemma 2.1]{B5} for a slightly general statement, but only for the case when ${\bf a}$ is decreasing and the case of ${\bf a}$ being increasing follows by applying the previous case to $-{\bf a}$.

   The above lemma allows us to deduce the following result:
\begin{theorem}
\label{thm4.2}
    Let $A, A'$ be two weighted mean matrix
  generated by $\{ \lambda_n \}^N_{n=1}$ and $\{ \lambda'_n \}^N_{n=1}$ respectively. Suppose
  that $\Lambda_n/\lambda_n \leq
  \Lambda'_n/\lambda'_n$ for all $n$. Then for fixed $p>1$, if \eqref{1.7} holds for
  $C=A$, we have $||A||_{p, p} \leq ||A'||_{p,
  p}$. Similarly, for fixed $p<0$, if $||A||^p_{p,p}$ is determined on an increasing
  sequence, we have $||A||^p_{p, p} \leq ||A'||^p_{p,
  p}$.
\end{theorem}
\begin{proof}
   Since the proofs are similar, we will only prove the $p>1$ case here. In this case as \eqref{1.7} holds for $C=A$, it follows from Lemma
   \ref{lem2} that $||A||_{p, p} \leq ||A'||_{p,
  p}$ as long as one can show that for any $k \leq n$,
\begin{equation*}
    \frac {\Lambda_k}{\Lambda_n} \leq \frac
    {\Lambda'_k}{\Lambda'_n}.
\end{equation*}
   By induction, it suffices to establish the above inequality for
   $k=n-1$ and one sees easily in this case the above inequality
   is equivalent to $\Lambda_n/\lambda_n \leq
  \Lambda'_n/\lambda'_n$ and this completes the proof.
\end{proof}

   We note here the above theorem can be regarded as in the
   spirit of Bennett's ``right is tight principle" (see page 409 of \cite{B5}) concerning the
   $l^p$ norms of summability matrices. According to the above theorem, we can interpret this
   principle for the weighted mean matrices as saying that for two
   given weighted mean matrices, the one with termwise larger diagonal
   entries has smaller norm, provided its norm is determined on
   decreasing sequences.

   As a concrete example of an application of the above theorem, we
   consider \eqref{8} for the cases $p>1, 1<\alpha<2$. As we mentioned earlier, if we
   assume \eqref{8} hold for those cases, then \eqref{1.7} holds for the
   corresponding matrix and in fact this is the case at least for
   $p \geq 2, 1< \alpha <2$ as \eqref{8} are known to hold for
   these cases. Now assume \eqref{1.7} does hold for the corresponding
   matrix for the cases $p>1, 1<\alpha<2$ of \eqref{8}, then in
   order to apply Theorem \ref{thm4.2} to establish \eqref{8}, we need to find
   a weighted mean matrix $A'$ (we may again focus on the finite matrices) whose $l^p$ norm is bounded by $\alpha p/(\alpha
   p-1)$. Now for the cases $1< \alpha <2$ of \eqref{8}, we consider the
following choice of
    the matrix $A'$ generated by $\{ \lambda'_n \}^N_{n=1}$,
    satisfying
\begin{equation*}
   \lambda'_1=1, \hspace{0.1in} \frac {\Lambda'_{n}}{\lambda'_{n}}=\frac {n+\alpha/2}{\alpha}, \hspace{0.1in} n \geq 2.
\end{equation*}
    Note that this defines the $\lambda'_n$'s uniquely and
    $\lambda'_n >0$ for all $n$. For a fixed $1<\alpha <2$, we now apply Theorem \ref{thm03} to conclude $||A'||_{p, p} \leq \alpha p/(\alpha
   p-1)$ for $p>1/(\alpha-1)^2$ by noting that it suffices to prove the case $n=1$ of \eqref{024} with $L=1/\alpha$ and this case follows when we bound
   $(1-1/(p\alpha))^{1-p}$ from below by $1-(1-p)/(p\alpha)+(1-1/p)/(2\alpha^2)$. It is also easy to check that for $n \geq 2$,
\begin{equation*}
   \frac {\sum^n_{k=1}k^{\alpha-1}}{n^{\alpha-1}} \leq \frac
   {n+\alpha/2}{\alpha }.
\end{equation*}
    One can similarly discuss the case $p<0, 1<\alpha<2$ using the following analogue of Theorem \ref{thm03}:
\begin{theorem}
\label{thm5.2}
    Let $p<0$ be fixed. Let $A$ be a weighted mean matrix generated by $\{ \lambda_n \}^N_{n=1}$. If for any integer $n \geq 1$, there exists a positive constant
    $L>0$ such that
\begin{equation*}
    \frac {\Lambda_{n+1}}{\lambda_{n+1}} \leq \frac
    {\Lambda_n}{\lambda_n}  \Big (1- \frac
    {L\lambda_n}{p\Lambda_n} \Big )^{1-p}+\frac {L}{p}~~,
\end{equation*}
    then $||A||^p_{p,p} \leq (p/(p-L))^p$.
\end{theorem}

   Apply the above theorem to $A'$ defined above, we see that $||A'||^p_{p, p} \leq (\alpha p/(\alpha p-1))^p$  and we then deduce immediately from Theorem \ref{thm4.2} the following
\begin{cor}
\label{thm4.3}
    Inequalities \eqref{8} hold for $p<0, 1 < \alpha <2$ for any
    increasing sequence ${\bf a}$.
\end{cor}

   Now, Corollary
   \ref{thm4.3} allows us to give another proof of the nontrivial cases $0< \alpha <1$ of Bennett's
   conjecture and in fact we shall prove a slightly general
   version by first establishing
\begin{theorem}
\label{thm4.4}
   Let $p<0$ be fixed and $N$ an integer and $A$ a weighted mean matrix generated by $\{ \lambda_n \}^{N}_{n=1}$. Suppose that the sequence
   $\{ \Lambda_n/\lambda_n \}^{\infty}_{n=1}$ is concave and that $\lim_{n
   \rightarrow +\infty}\Lambda_n/(n\lambda_n)=L$. If we have
\begin{equation}
\label{4.0}
   e^{\lambda_1/\lambda_2}(1-L) <1,
\end{equation}
  then $||A||^p_{p,p}$ is determined on an increasing
  sequence.
\end{theorem}
\begin{proof}
   As $\{ \Lambda_n/\lambda_n \}^{\infty}_{n=1}$ is concave and that $\lim_{n
   \rightarrow +\infty}\Lambda_n/(n\lambda_n)=L$, a result of
   Bennett \cite[Lemma 2]{Be1} implies that
   $L \leq \Lambda_{n+1}/\lambda_{n+1}-\Lambda_{n}/\lambda_n \leq \Lambda_{2}/\lambda_{2}-\Lambda_{1}/\lambda_1=\lambda_1/\lambda_2$. It
   follows from Theorem \ref{thm1'} that $||A||^p_{p,p} \leq (p/(p-\lambda_1/\lambda_2))^p$
   for $p<0$. Thus inequality \eqref{2.4} holds with $U_{p, N}=(p/(p-\lambda_1/\lambda_2))^p$. As $\lim_{p \rightarrow
   -\infty}(p/(p-\lambda_1/\lambda_2))^p=e^{\lambda_1/\lambda_2}$
   and $(p/(p-\lambda_1/\lambda_2))^p$ is a decreasing function
   of $p<0$, we see that inequality \eqref{1.8} holds with
   $U_{p, N}=(p/(p-\lambda_1/\lambda_2))^p$ by \eqref{4.0}. Now
   our assertion follows from Theorem \ref{thm4}.
\end{proof}

   We now apply the above theorem to $\lambda_n=n^{\alpha}$ for $0 <\alpha<1$, in which case \eqref{4.0} is equivalent to
\begin{equation*}
    1+\frac 1{\alpha} > e^{1/2^{\alpha}}.
\end{equation*}
    As $e^{1/2^{\alpha}}<e$ when $0< \alpha <1$, it follows that
    the above inequality holds for $\alpha < 1/(e-1) \approx 0.58$.
    Thus we may assume that $1/2 \leq \alpha <1$ and in this case
    $e^{1/2^{\alpha}}<e^{1/\sqrt{2}}$ and by repeating the above
    argument, we see that we may further assume that $0.8 \leq \alpha
    <1$ but then the above inequality holds since $e^{1/2^{0.8}}
    <2$. Therefore, combined with Corollary \ref{thm4.3}, we see that
    inequalities \eqref{8} hold for $p<0, 1 < \alpha <2$ and for
    the other positive $\alpha$'s, we can apply Theorem \ref{thm1'} to
    conclude that inequalities \eqref{8} hold as well and we
    summarize our result in the following
\begin{cor}
    Inequalities \eqref{8} hold for $p<0, \alpha >0$.
\end{cor}

    We note here that the above corollary implies the nontrivial cases $0< \alpha <1$ of Bennett's
   conjecture, which one obtains by taking $p \rightarrow -\infty$
   of the corresponding cases of \eqref{8}.

\section{A generalization of a result of Bennett }
\label{sec 4} \setcounter{equation}{0}
  As we mentioned in the introduction, the validity of \eqref{1.7} will allow us to deduce the cases $1/p <
\alpha \leq 1$ of inequalities \eqref{7}. In this section, we
shall generalize a result of Bennett which in turn implies these
cases. We shall assume all the infinite sums converge and we start
by noting the following result of Bliss \cite{Bl}:
\begin{theorem}
\label{thm6}
  Let $r > p > 1$ and let $\alpha$ be a real number satisfying $(\alpha+1)p>1$. Let $f(x)$ be a non-negative measurable function on $[0, +\infty)$ such that $f \in L^p(0, +\infty)$. Then the integral $\int^x_0f(t)t^{\alpha}dt$ is finite for every $x$ and
\begin{equation*}
  \int^{\infty}_{0}\Big ( \int^x_0f(t)t^{\alpha}dt \Big )^r\frac {dx}{x^{(\alpha+1)r-s}} \leq K_{r,s, \alpha}\Big ( \int^{\infty}_0f(x)^pdx \Big )^{r/p},
\end{equation*}
  where
\begin{equation*}
  s=r/p-1, \hspace{0.1in} K_{r,s,\alpha}=\frac {1}{(r-s-1)(1+\alpha q)^{r-s}}\Big ( \frac {s\Gamma(r/s)}{\Gamma(1/s)\Gamma((r-1)/s)} \Big )^s.
\end{equation*}
\end{theorem}

   We note here Bliss only proved the case $\alpha=0$ in \cite{Bl}
   but the general case can be obtained by some changes of variables.
   Based on the above result, we now prove the following
\begin{theorem}
\label{thm7}
  Let $r > s > 1$ and $s/r < \alpha \leq 1$. Let ${\bf u}, {\bf v}, {\bf a}$ be sequences with positive entries.
  Let $V_n=\sum^n_{k=1}v_k$ for $n \geq 1$ and $V_0=0$. If for $m
  \geq 1$,
\begin{equation*}
  \sum^m_{n=1}u_nV^{\alpha r}_n \leq V^s_m.
\end{equation*}
   Then
\begin{equation*}
  \sum^\infty_{n=1}u_n\Big ( \sum^n_{k=1}(V^{\alpha}_k-V^{\alpha}_{k-1})a_k \Big )^{r} \leq
  s\alpha^rK_{r,s-1, \alpha-1}\Big ( \sum^\infty_{n=1}v_na^{r/s}_n \Big )^{s}.
\end{equation*}
\end{theorem}
\begin{proof}
   The proof is almost identical to the proof of Theorem 2 in
   \cite{B3}, taking account into Theorem \ref{thm6}, as long as
   one can show (see also the proof of Theorem 1 in \cite{B1}) that for $1 \leq i<j, a_i < a_j$,
\begin{equation*}
  \frac {(V^{\alpha}_i-V^{\alpha}_{i-1})a_i+(V^{\alpha}_j-V^{\alpha}_{j-1})a_j}
  {(V^{\alpha}_i-V^{\alpha}_{i-1})+(V^{\alpha}_j-V^{\alpha}_{j-1})}
  \leq \frac {v_ia_i+v_ja_j}
  {v_i+v_j}.
\end{equation*}
   The above inequality follows from Lemma \ref{lem2} (note that
   $a_j>a_i$ here) provided that
\begin{equation*}
  \frac {V^{\alpha}_j-V^{\alpha}_{j-1}}{V_j-V_{j-1}}
    \leq \frac {V^{\alpha}_i-V^{\alpha}_{i-1}}{V_i-V_{i-1}}.
\end{equation*}
  The above inequality holds by the mean value theorem, since
  the right-hand side is no less than $\alpha V^{\alpha-1}_i$ and
  the left-hand side is no greater than $\alpha
  V^{\alpha-1}_{j-1}$ and this completes the proof.
\end{proof}

   We now take $u_n=(V^{s}_n-V^{s}_{n-1})/V^{\alpha r}_n$ in the
   above theorem and make a change of variables $a_n \rightarrow a^{s/r}_n$ and let $r \rightarrow +\infty$ to deduce that
\begin{cor}
\label{cor7}
  Let $s > 1$ and $0 < \alpha \leq 1$. Let ${\bf v}, {\bf a}$ be sequences with positive entries.
  Let $V_n=\sum^n_{k=1}v_k$ for $n \geq 1$ and $V_0=0$.
   Then
\begin{equation*}
  \sum^\infty_{n=1}(V^{s}_n-V^{s}_{n-1})\Big ( \prod^n_{k=1}a^{V^{\alpha}_k-V^{\alpha}_{k-1}}_k \Big )^{s/V^{\alpha}_n} \leq
  \frac {e^{-(\alpha-1)s/\alpha}}{\alpha^{1-s}}\frac {s}{s-1}\Big ( \frac {s-1}{\Gamma(1/(s-1))} \Big )^{s-1}\Big ( \sum^\infty_{n=1}v_na_n \Big )^{s}.
\end{equation*}
\end{cor}
  Note that we will get back Carleman-type inequalities on letting
  $s \rightarrow 1^+$ in the above corollary. We can also take $v_n=1$ and $u_n=(n^{s}-(n-1)^{s})/n^{\alpha
  r}$ in Theorem \ref{thm6} to deduce that
\begin{cor}
\label{cor8}
  Let $r>s > 1$ and $s/r < \alpha \leq 1$. Let ${\bf a}$ be sequences with positive entries.
   Then
\begin{equation*}
  \sum^\infty_{n=1}(n^{s}-(n-1)^{s})\Big ( \frac {1}{n^{\alpha}} \sum^n_{k=1}(k^{\alpha}-(k-1)^{\alpha})a_k \Big )^{r} \leq
  s\alpha^rK_{r,s-1, \alpha-1}\Big ( \sum^\infty_{n=1}v_na^{r/s}_n \Big
 )^{s}.
\end{equation*}
\end{cor}
  Note that we get back the cases $1/p < \alpha \leq 1$ of
  \eqref{7} on setting $r=p$ and letting
  $s \rightarrow 1^+$ in the above corollary.

\section{Schur's Test and Some Generalizations of inequalities \eqref{7} and \eqref{8}}
\label{sec 8} \setcounter{equation}{0}
   In this section we first state a discrete version of Schur's test concerning the norms of linear
   operators:
\begin{lemma}
\label{lem2.1} Let $p>1$ be fixed and let $A=(\alpha_{j,i})_{1 \leq i, j \leq N}$ be a matrix with non-negative entries. If there exist positive numbers $U_1, U_2$ and two positive sequences ${\bf c}=(c_i), 1 \leq i \leq N; {\bf d}=(d_i), 1 \leq i \leq N$, such that
\begin{eqnarray}
\label{8.1}
  \sum^{N}_{i=1} \alpha_{j,i}c^{1/p}_i  &\leq & U_1d^{1/p}_j,  \hspace{0.1in}  1 \leq j \leq N; \\
\label{8.2}
 \sum^{N}_{j=1}\alpha_{j,i}d^{1/q}_j &\leq & U_2c^{1/q}_i,  \hspace{0.1in}  1 \leq i \leq N.
\end{eqnarray}
   Then
\begin{equation*}
  ||A||_{p,p} \leq U^{1/q}_1U^{1/p}_2.
\end{equation*}
\end{lemma}


   We now point out that Schur's test is equivalent to the approaches mentioned in \cite{G5} in determining the operator norms of weighted mean matrices.
   It suffices to show that it is equivalent to the approach of Kaluza and
Szeg\"o. To see this, note that our goal in general is to find
some (smallest possible) constant $U_{p, N}$ so that for a
weighted mean matrix generated by $\{ \lambda_n \}^{N}_{n=1}$ (we
may assume $\lambda_n>0$ for all $n$), inequality \eqref{3.1}
holds for any integer $N \geq 1$ and any ${\bf a} \in l^p$.
   We now apply Lemma \ref{lem2.1} with $\alpha_{j,i}=\lambda_i/\Lambda_j$ for $i \leq j$ and $\alpha_{j,i}=0$ for $i > j$ with
\begin{equation*}
   c_i=\Big ( \frac {w_i}{\lambda_i }\Big )^p, \hspace{0.1in} d_j=\Big ( \frac {\sum^{j}_{k=1}w_k }{\Lambda_j} \Big )^p,
\end{equation*}
  where the auxiliary sequence $\{ w_n \}^{\infty}_{n=1}$ is of positive terms and to be determined later.
  The choice of the $c_i$'s and $d_j$'s is to make inequality \eqref{8.1} satisfied with $U_1=1$ (it becomes an identity) and inequality \eqref{8.2} becomes
\begin{equation}
\label{8.3}
  \sum_{j=i}^{N}\frac
  {\lambda_i}{\Lambda^p_j}\Big(\sum_{k=1}^jw_k\Big)^{p-1} \leq U_2 \Big ( \frac {w_i}{\lambda_i }\Big )^{p-1}.
\end{equation}
  Suppose now one can find for each $p>1$ a positive constant $U_2$,
  a sequence ${\bf w}$ of positive terms with $w_n^{p-1}/\lambda^p_n$ decreasing to $0$, such
  that for any integer $n \geq 1$,
\begin{equation*}
 (w_1+\cdots+w_n)^{p-1}< U_2\Lambda^p_n( \frac {w_n^{p-1}}{\lambda^p_n}-\frac {w_{n+1}^{p-1}}{\lambda^p_{n+1}} ),
\end{equation*}
  then inequality \eqref{8.3} will follow from this and this is exactly the starting point of Kaluza and
Szeg\"o's approach.


   In what follows, we will give an account of Bennett's proof of \eqref{7} in the form of Schur's test.
   First we consider the case $\alpha > 1/p$ of \eqref{7} and we can replace the infinite sums by  finite sums from $1$ to $N$
   with $N \geq 1$  here and we note the following estimation (\cite[(99)]{Be1}):
\begin{equation}
\label{2.5}
 \sum^{N}_{j=i}\frac {\int^{i}_{i-1}x^{\alpha-1/p}dx }{j^{\alpha+1/q}} \leq  \frac {1}{\alpha -1/p}.
\end{equation}
   We now apply Lemma  \ref{lem2.1} with
$\alpha_{j,i}=\alpha \Big ( \int^i_{i-1}x^{\alpha-1/p}dx \Big
)^{1/p}\Big( \int^i_{i-1}x^{\alpha-1/p-1}dx \Big
)^{1/q}/j^{\alpha}$ for $i \leq j$ and $a_{j,i}=0$ otherwise and
$c_i= (\int^i_{i-1}x^{\alpha-1/p-1}dx/
\int^{i}_{i-1}x^{\alpha-1/p}dx), d_j=1/j$, $U_1=U_2=(\alpha
p)/(\alpha p-1)$
 to see that in this case inequality \eqref{8.1} becomes an identity and inequality \eqref{8.2} becomes exactly \eqref{2.5}. From this we deduce the following inequality for $p>1, \alpha > 1/p$ and any ${\bf a} \in l^p$,
\begin{equation*}
 \sum^{\infty}_{n=1}\Big{|}\frac
1{n^{\alpha}}\sum^n_{i=1}\alpha \Big ( \int^i_{i-1}x^{\alpha-1/p}dx \Big )^{1/p}\Big ( \int^i_{i-1}x^{\alpha-1-1/p}dx \Big )^{1/q}a_i\Big{|}^p
\leq  \Big( \frac {\alpha p}{\alpha p-1} \Big )^p\sum^{\infty}_{n=1}|a_n|^p.
\end{equation*}
 from which one deduces the corresponding cases of \eqref{7} easily.


  We note here in Bennett's proof of \eqref{7} given above, a key ingredient is inequality \eqref{2.5}.
  We point out here that when $1 \leq \alpha \leq 1+1/p$, a better estimation exists, namely,
\begin{equation}
\label{4.10}
 \sum^{N}_{j=i}\frac {\alpha \Big(i-\frac {1}{2} \Big )^{\alpha-1+1/q}}{j^{\alpha+1/q}} \leq  \frac {\alpha p}{\alpha p-1}.
\end{equation}
  Inequality \eqref{4.10} can be easily deduced from the following inequality for all integers $i \geq 1$ and $1 \leq \alpha \leq 1+1/p$,
\begin{equation*}
  i^{-\alpha-1/q} \leq  \frac {1}{\alpha-1/p} \Big ( \Big ( i-1/2 \Big )^{1-\alpha-1/q}-\Big ( i+1/2 \Big )^{1-\alpha-1/q} \Big )=\int^{i+1/2}_{i-1/2}x^{-\alpha-1/q}dx.
\end{equation*}
   The above inequality follows from the well-known Hadamard's inequality (with $h(x)=x^{-\alpha-1/q}, a=i-1/2, b=i+1/2$ below), which asserts that for a continuous convex function $h(x)$ on $[a, b]$,
\begin{equation*}
   h(\frac {a+b}2) \leq \frac {1}{b-a}\int^b_a h(x)dx \leq \frac
   {h(a)+h(b)}{2}.
\end{equation*}
   The above inequality also allows us to see easily that inequality
   \eqref{4.10} improves upon \eqref{2.5} for $1 \leq \alpha \leq
   1+1/p$.

  Now, inequality \eqref{4.10} allows us to establish the following
\begin{theorem}
\label{thm4.1}
    Let $p>1$ be fixed, then the following inequality holds for $1 \leq \alpha \leq 1+1/p$ and any ${\bf a} \in l^p$,
 \begin{equation*}
 \sum^{\infty}_{n=1}\Big{|}\frac
1{n^{\alpha}}\sum^n_{i=1}\alpha \Big(i-\frac {1}{2} \Big )^{\frac {1}{p}(\alpha-\frac {1}{p})}\Big ( \int^i_{i-1}x^{\alpha-1-1/p}dx \Big )^{1/q}a_i\Big{|}^p
\leq  \Big( \frac {\alpha p}{\alpha p-1} \Big )^p\sum^{\infty}_{n=1}|a_n|^p.
\end{equation*}
\end{theorem}
\begin{proof}
    We can replace the infinite sums by finite sums from $1$ to $N$ with $N \geq 1$  here
    and we apply Lemma  \ref{lem2.1} here
    with $\alpha_{j,i}=\alpha \Big(i-\frac {1}{2} \Big )^{\frac {1}{p}(\alpha-\frac {1}{p})}
    \Big ( \int^i_{i-1}x^{\alpha-1-1/p}dx \Big )^{1/q}/j^{\alpha}$
    for $i \leq j$ and $0$ otherwise and $c_i=\Big(i-\frac {1}{2} \Big )^{-(\alpha-\frac {1}{p})}\Big ( \int^i_{i-1}x^{\alpha-1-1/p}dx \Big ),
    d_j=j^{-1}$ to see that estimations \eqref{8.1}-\eqref{8.2} hold by \eqref{4.10}
    with $U_1=U_2=\alpha p/(\alpha p-1)$ and this completes the proof.
\end{proof}

  To deduce interesting corollaries from Theorem \ref{thm4.1}, we note the following lemma:
\begin{lemma}[{\cite[Lemma 2.1]{alz1.5}}]
\label{lem4}
   Let $a > 0, b > 0$ and $r$ be real numbers with $a \neq b$, and
   let
\begin{eqnarray*}
   L_r (a, b) &=& \genfrac(){1pt}{}{a^r - b^r}{r(a-b)}^{1/(r-1)}  \hspace{0.2in} (r \neq  0,
   1), \\
   L_0(a,b) &=& \frac {a - b}{ \log a-\log
  b}, \\
   L_1(a,b) &=& \frac 1{e}\genfrac(){1pt}{}{a^a}{b^b}^{1/(a-b)}.
\end{eqnarray*}
   The function $r \mapsto L_r(a,b)$ is strictly increasing on ${\mathbb R}$.
\end{lemma}
   It readily follows from the above lemma that for $1 \leq \alpha \leq 1+1/p$, we have
\begin{eqnarray*}
  && i^{\alpha}-(i-1)^{\alpha} = \alpha L^{\alpha-1}_{\alpha}(i, i-1)  \leq \alpha L^{\alpha-1}_{2}(i, i-1) = \alpha \Big ( i-1/2 \Big )^{\alpha-1} \\
 & \leq & \alpha \Big(i-1/2 \Big )^{\frac {1}{p}(\alpha-\frac {1}{p})}\Big ( \int^i_{i-1}x^{\alpha-1-1/p}dx \Big )^{1/q}=\alpha L^{\frac {1}{p}(\alpha-\frac {1}{p})}_{2}(i, i-1) \cdot L^{\frac {1}{q}(\alpha-1-\frac {1}{p})}_{\alpha-\frac {1}{p}}(i, i-1).
\end{eqnarray*}
  It follows from this that Theorem \ref{thm4.1} not only implies the corresponding cases of \eqref{7} but also the following stronger result:
\begin{cor}
\label{cor4.1}
    Let $p>1$ be fixed, then the following inequality holds for $1 \leq \alpha \leq 1+1/p$ and any ${\bf a} \in l^p$,
 \begin{equation*}
 \sum^{\infty}_{n=1}\Big{|}\frac
1{n^{\alpha}}\sum^n_{i=1}\alpha \Big(i-\frac {1}{2} \Big )^{\alpha-1}a_i\Big{|}^p
\leq  \Big( \frac {\alpha p}{\alpha p-1} \Big )^p\sum^{\infty}_{n=1}|a_n|^p.
\end{equation*}
\end{cor}

  As an interesting consequence of Corollary \ref{cor4.1}, we note for the case $p=2$ we have $2\alpha -1 \leq 2$ for $\alpha \leq 3/2$
  so that Corollary \ref{cor4.1} implies the following inequality for ${\bf a} \in l^2$ and $1 \leq \alpha \leq 3/2$:
\begin{equation}
\label{4.4}
  \sum^{N}_{n=1}\Big | \sum^{n}_{i=1}\frac {\alpha L^{\alpha-1}_{2\alpha-1}(i, i-1)}{n^{\alpha}}a_i\Big |^2 \leq \frac {\alpha^2}{(\alpha-1/2)^2}\sum^{N}_{i=1}|a_i|^2.
\end{equation}

   We now apply the duality principle \cite[Lemma 2]{M} to deduce from \eqref{4.4} the following inequality for ${\bf a} \in l^2, a_i \geq 0$ and $1 \leq \alpha \leq 3/2$:
\begin{equation*}
   \sum^{N}_{i,j=1}\frac {\alpha^2 \min ( i^{2\alpha-1}, j^{2\alpha-1}) }{(2\alpha-1)i^{\alpha}j^{\alpha}}a_ia_j= \sum^{N}_{n=1}\Big ( \sum^{N}_{i=n}\frac {\alpha L^{\alpha-1}_{2\alpha-1}(n, n-1)}{i^{\alpha}}a_i\Big )^2 \leq \frac {\alpha^2}{(\alpha-1/2)^2}\sum^{N}_{i=1}a^2_i.
\end{equation*}

   We note here the case $\alpha=1$ above gives back  a result of Schur in \cite {Schur}, who showed that for ${\bf x}, {\bf y} \in l^2$,
\begin{equation*}
   \sum^{\infty}_{i,j=1}\frac {x_iy_j}{\max (i,j)} \leq 4 || {\bf x}||_2 || {\bf y}||_2.
\end{equation*}
   By the duality principle, the above inequality is equivalent to Hardy's inequality \eqref{eq:1} for the case $p=2$, even though this was not mentioned in \cite{Schur} (this is actually prior to Hardy's discovery of \eqref{eq:1}).

   Our discussions above allow us to regard the cases of $\alpha \geq 1$ of inequalities \eqref{7} and \eqref{8}
   as special cases of a family of inequalities. Namely, it is interesting to determine the best constant $U=U(\alpha, \beta, p)$ so that
   the following inequality holds for all ${\bf a} \in l^p$ ($p>1, \beta \geq \alpha \geq 1$):
\begin{equation}
\label{4.5}
 \sum^{\infty}_{n=1}\Big{|}\frac
1{\sum^n_{k=1}L^{\alpha-1}_{\beta}(k,
k-1)}\sum^n_{i=1}L^{\alpha-1}_{\beta}(i, i-1)a_i \Big{|}^p \leq
U\sum^{\infty}_{n=1}|a_n|^p.
\end{equation}
  Note that the case of $\beta=\alpha$ above corresponds to inequality \eqref{7} and the case of
$\beta \rightarrow +\infty$ above corresponds to inequality
\eqref{8} by Lemma \ref{lem4}. In both cases, we expect $U=(\alpha
p /(\alpha p-1) )^p$ (of course this is known except for some
cases of \eqref{8} when $1<p<2, 1<\alpha<2$. Thanks to Corollary
\ref{cor4.1} and Lemma \ref{lem4}, we also know that inequality
\eqref{4.5} holds with $U=(\alpha p /(\alpha p-1) )^p$ for $p>1, 1
\leq \alpha \leq 1+1/p, \alpha \leq \beta \leq 2$.

\section*{Acknowledgement}
  The author is supported by a research fellowship from an Academic
Research Fund Tier 1 grant at Nanyang Technological University for
this work.


\end{document}